# On prediction errors in regression models with nonstationary regressors

## Ching-Kang Ing[1] and Chor-Yiu Sin[2]

*Academia Sinica and Xiamen University*

**Abstract:** In this article asymptotic expressions for the final prediction error (FPE) and the accumulated prediction error (APE) of the least squares predictor are obtained in regression models with nonstationary regressors. It is shown that the term of order $1/n$ in FPE and the term of order $\log n$ in APE share the same constant, where $n$ is the sample size. Since the model includes the random walk model as a special case, these asymptotic expressions extend some of the results in Wei (1987) and Ing (2001). In addition, we also show that while the FPE of the least squares predictor is not affected by the contemporary correlation between the innovations in input and output variables, the mean squared error of the least squares estimate does vary with this correlation.

## 1. Introduction

Consider a simple regression model

$$(1.1) \qquad y_t = \beta x_{t-1} + \varepsilon_t,$$

where $\beta$ is an unknown constant, $\varepsilon_t$'s are (unobservable) independent random disturbances with zero means and a common variance $\sigma^2$, and $x_t$ is an unit root process satisfying

$$(1.2) \qquad x_t = x_{t-1} + \eta_t,$$

with $x_0 = 0$, $\eta_t = \sum_{j=0}^{t-1} c_j \omega_{t-j}$, $\sum_{j=0}^{\infty} |c_j| < \infty$, $\sum_{j=0}^{\infty} c_j \neq 0$, and $\omega_t$ being independent random noises with zero means and a common variance $\sigma_\omega^2$. We also assume that $\varepsilon_t$ is independent of $\{\omega_j, j \leq t-1\}$. Note that if $\beta = 1$, $c_0 = 1$, $c_j = 0$ if $j > 0$, and $\varepsilon_t = \omega_t$, then (1.1) becomes the well-known random walk model (see, for instance, Chan and Wei [4]). Having observed $(y_{i+1}, x_i), i = 1, \ldots, n-1$, $\beta$ can be estimated by least squares

$$(1.3) \qquad \hat{\beta}_n = \frac{\sum_{i=1}^{n-1} x_i y_{i+1}}{\sum_{i=1}^{n-1} x_i^2}.$$

If $x_n$ also becomes available, then it is natural to predict $y_{n+1}$ using the least squares predictor,

$$(1.4) \qquad \hat{y}_{n+1} = x_n \hat{\beta}_n.$$

---

[1]Academia Sinica, Taipei, Taiwan, R.O.C., e-mail: cking@stat.sinica.edu.tw
[2]Wang Yanan Institute for Studies in Economics, Xiamen University, Fujian, e-mail: cysinhkbu@gmail.com

To assess the performances of the least squares predictor, we consider the final prediction error (FPE, Akaike [1])

$$(1.5) \qquad E\left\{(y_{n+1} - \hat{y}_{n+1})^2\right\} = \sigma^2 + E\left\{x_n^2(\hat{\beta}_n - \beta)^2\right\},$$

and the accumulated prediction error (APE, Rissanen [14])

$$\sum_{i=2}^{n}(y_i - \hat{y}_i)^2 = \sum_{i=2}^{n}\left\{\varepsilon_i - x_{i-1}(\hat{\beta}_{i-1} - \beta)\right\}^2$$

$$(1.6) \qquad = \sum_{i=2}^{n}\varepsilon_i^2 + \sum_{i=2}^{n}x_{i-1}^2(\hat{\beta}_{i-1} - \beta)^2(1 + o(1)) \text{ a.s.},$$

where the second equality of (1.6) is ensured by Chow [5]. It is straightforward to see that the terms in (1.5) and (1.6),

$$(1.7) \qquad \sum_{i=2}^{n} x_{i-1}^2(\hat{\beta}_{i-1} - \beta)^2 = \sum_{i=2}^{n}\left\{\frac{x_i^2(\sum_{j=1}^{i-1} x_j \varepsilon_{j+1})^2}{(\sum_{j=1}^{i-1} x_j^2)^2}\right\},$$

and

$$(1.8) \qquad n x_n^2(\hat{\beta}_n - \beta)^2 = \left\{\frac{(\frac{1}{\sqrt{n}} x_n)(\frac{1}{n}\sum_{i=1}^{n-1} x_i \varepsilon_{i+1})}{\frac{1}{n^2}\sum_{i=1}^{n-1} x_i^2}\right\}^2.$$

When $\{y_t\}$ is a random walk model mentioned above, Wei ([15], Theorem 4) showed that the rhs of (1.7) equals $2\sigma_\omega^2 \log n + o(\log n)$ a.s. By imposing further assumptions on the distribution of $\omega_t$, Ing ([9], Corollary 1) subsequently obtained the limiting value of the expectation on the rhs of (1.8), which is $2\sigma_\omega^2$. This article extends these two results to models (1.1) and (1.2), which provides a deeper understanding of the least squares predictor (estimate) in situations where Fisher's information, $\sum_{j=1}^{n-1} x_j^2$, grows at a rate much faster than $n$, and the innovations in input and output variables come from different sources. The rest of the paper is organized as follows. Section 2 derives the asymptotic expressions for the rhs of (1.7). In Section 3, sufficient conditions are given to ensure that the expectation on the rhs of (1.8) is bounded by some finite positive constant for all sufficiently large $n$. We then apply this moment property and the results obtained in Section 2 to show that

$$(1.9) \qquad \lim_{n \to \infty} E\{n x_n^2(\hat{\beta}_n - \beta)^2\} = 2\sigma^2.$$

Some discussions related to (1.9) are given at the end of Section 3. In particular, it is shown that while the FPE of the least squares predictor is not affected by the contemporary correlation between $\varepsilon_t$ and $\omega_t$, the mean squared error of the least squares estimate does vary with this correlation. In addition, we also show that the squares of the normalized estimate, $n(\hat{\beta}_n - \beta)$, and the normalized regressor, $x_n/\sqrt{n}$, are not asymptotically uncorrelated.

## 2. An asymptotic expression for the APE

To prove the main result of this section, two auxiliary lemmas are required. They are also of independent interests.

**Lemma 1.** *Assume the $\{\omega_t\}$ in Section 1 satisfy $\sup_{-\infty<t<\infty} E|\omega_t|^\alpha < \infty$ for some $\alpha > 2$. Let $z_t = \sum_{j=0}^{t-1} d_j \omega_{t-j}$, where $|d_j| \leq Cj^{-1}$ for some $C > 0$ and all $j \geq 1$. Then, with $\gamma_t = \sigma_\omega^2 \sum_{j=0}^{t-1} d_j^2$,*

$$\frac{1}{n}\sum_{t=1}^{n}(z_t^2 - \gamma_t) = o(1) \text{ a.s.} \tag{2.1}$$

*Proof.* Straightforward calculations yield that

$$z_t^2 - \gamma_t = \sum_{l=1}^{t} d_{t-l}^2(\omega_l^2 - \sigma_\omega^2) + 2\sum_{l_2=2}^{t}\sum_{l_1=1}^{l_2-1} d_{t-l_1}d_{t-l_2}\omega_{l_1}\omega_{l_2}. \tag{2.2}$$

By (2.2) and changing the order of summations,

$$\sum_{t=n_1}^{n_2} \frac{z_t^2 - \gamma_t}{t} = \sum_{l=1}^{n_1}\left(\sum_{t=n_1}^{n_2} \frac{d_{t-l}^2}{t}\right)\eta_l^* + \sum_{l=n_1+1}^{n_2}\left(\sum_{t=l}^{n_2} \frac{d_{t-l}^2}{t}\right)\eta_l^*$$
$$+ 2\sum_{l_2=2}^{n_1}\left\{\sum_{l_1=1}^{l_2-1}\left(\sum_{t=n_1}^{n_2} \frac{d_{t-l_1}d_{t-l_2}}{t}\right)\omega_{l_1}\right\}\omega_{l_2}$$
$$+ 2\sum_{l_2=n_1+1}^{n_2}\left\{\sum_{l_1=1}^{l_2-1}\left(\sum_{t=l_2}^{n_2} \frac{d_{t-l_1}d_{t-l_2}}{t}\right)\omega_{l_1}\right\}\omega_{l_2}$$
$$\equiv (1) + (2) + (3) + (4),$$

where $\eta_t^* = \omega_t^2 - \sigma_\omega^2$. In the following, we shall show that for some $\alpha_k > 1$, there are $C_k > 0, \xi_{1,k} > 1$, and $\xi_{2,k} > 1$ independent of $n_1$ and $n_2$ such that

$$E|(k)|^{\alpha_k} \leq C_k \Big(\sum_{t=n_1}^{n_2} \frac{1}{t^{\xi_{1,k}}}\Big)^{\xi_{2,k}}, \tag{2.3}$$

where $k = 1,\ldots,4$. (2.3) and Móricz (1976) imply that for some $\alpha > 1$, there are $C^* > 0, \xi_1 > 1$, and $\xi_2 > 1$ independent of $n_1$ and $n_2$ such that

$$E \max_{n_1 \leq l \leq n_2} \Big|\sum_{t=n_1}^{l} \frac{z_t^2 - \gamma_t}{t}\Big|^\alpha \leq C^*\Big(\sum_{t=n_1}^{n_2} \frac{1}{t^{\xi_1}}\Big)^{\xi_2}. \tag{2.4}$$

As a result, (2.1) follows from (2.4) and Kronecker's lemma.

Let $\alpha_1 = \min\{\alpha/2, 2\}$. Then,

$$E|(1)|^{\alpha_1} \leq C_{1,1} E\Big\{\sum_{l=1}^{n_1}\Big(\sum_{t=n_1}^{n_2} \frac{d_{t-l}^2}{t}\Big)^2 \eta_l^{*2}\Big\}^{\alpha_1/2}$$
$$\leq C_{1,1} \sum_{t_1=n_1}^{n_2}\sum_{t_2=n_1}^{n_2} \frac{1}{t_1^{\alpha_1/2} t_2^{\alpha_1/2}} \sum_{l=1}^{n_1} |d_{t_1-l}d_{t_2-l}|^{\alpha_1} E|\eta_l^*|^{\alpha_1}$$
$$\leq C_{1,2}\left(\sum_{t=n_1}^{n_2} \frac{1}{t^{\alpha_1}} + \sum_{t_1=n_1}^{n_2-1} \frac{1}{t_1^{\alpha_1/2}} \sum_{t_2=t_1+1}^{n_2} \frac{1}{t_2^{\alpha_1/2}}(t_2-t_1)^{-\alpha_1}\right) \tag{2.5}$$
$$\leq C_{1,3}\left(\sum_{t=n_1}^{n_2} \frac{1}{t^{\alpha_1}}\right) \leq C_{1,3}\left(\sum_{t=n_1}^{n_2} \frac{1}{t^{\xi_{1,1}}}\right)^{\xi_{2,1}},$$

where $C_{1,i}, i = 1, 2, 3$ are some positive constant independent of $n_1$ and $n_2$, $1 < \xi_{1,1} < \alpha_1$, $\xi_{2,1} = \alpha_1/\xi_{1,1}$, first inequality follows from Burkholder's inequality, second one follows from the fact that $0 < \alpha_1/2 \leq 1$ and changing the order of summations, third one is ensured by $\sup_t E|\omega_t|^\alpha < \infty$ and $|d_j| \leq Cj^{-1}$, which implies for all $n_1 \leq t_1, t_2 \leq n_2$, $\sum_{l=1}^{n_1} |d_{t_1-l}d_{t_2-l}|^{\alpha_1} \leq C_{1,4}|t_1 - t_2|^{-\alpha_1}$, for some $C_{1,4} > 0$. As a result, (2.3) holds for $k = 1$. The proof of (2.3) for the case of $k = 2$ is similar. The details are thus omitted. To show (2.3) for the case $k = 3$, let $\alpha_3 = \alpha$. Then, by Minkowski's inequality and using Wei (1987, Lemma 2) twice, one obtains

$$(2.6) \quad \begin{aligned} E|(3)|^{\alpha_3} &\leq C_{3,1} E \left| \sum_{l_2=2}^{n_1} \left\{ \sum_{l_1=1}^{l_2-1} \left( \sum_{t=n_1}^{n_2} \frac{d_{t-l_1}d_{t-l_2}}{t} \right) \omega_{l_1} \right\} \omega_{l_2} \right|^{\alpha_3} \\ &\leq C_{3,2} \left( \sum_{l_2=2}^{n_1} \sum_{l_1=1}^{l_2-1} (\sum_{t=n_1}^{n_2} \frac{d_{t-l_1}d_{t-l_2}}{t})^2 \right)^{\alpha_3/2}, \end{aligned}$$

where $C_{3,i}, i = 1, 2$ are some positive constants independent of $n_1$ and $n_2$. Observe that for $n_1 \leq t_1 < t_2 \leq n_2$ and any $1 \leq M_1 \leq M_2 \leq n_1$, $\sum_{l=M_1}^{M_2} |d_{t_1-l}d_{t_2-l}| \leq C_{3,3}(\log t_2 - \log t_1)/(t_2 - t_1)$, where $C_{3,3} > 0$ is independent of $M_1$ and $M_2$. Using this fact and changing the order of summations, it follows that the rhs of (2.6) is bounded by $C_{3,4}(\sum_{t=n_1}^{n_2} t^{-2})^{\alpha_3/2}$, where $C_{3,4}$ is a positive constant independent of $n_1$ and $n_2$. Hence, (2.3) holds for $k = 3$. The proof of (2.3) for the case $k = 4$ is similar to that of $k = 3$. Therefore, we skip the details. □

**Remark 1.** If in Lemma 1 $z_t = \sum_{j=0}^{\infty} d_j \omega_{t-j}$ with $|d_j| \leq Cj^{-1}, j \geq 1$, then the same argument also yields (2.1) but with $\gamma_t$ replaced by $\gamma^* = \sigma_\omega^2 \sum_{j=0}^{\infty} d_j^2$. For a related result, Brockwell and Davis (1987, Proposition 7.3.5), assuming that $\omega_j$'s are i.i.d. with finite second moment and $d_j$'s satisfy $\sum_{j=0}^{\infty} |d_j| < \infty$ and $\sum_{j=0}^{\infty} d_j^2 j < \infty$, obtained $(n^{-1} \sum_{t=1}^{n} z_t^2) - \gamma^* = o_p(1)$. While the moment restriction of their result is slightly weaker than that of Lemma 1, the identically distributed assumption can be dropped in Lemma 1. In addition, the assumption on $d_j$ in Lemma 1 seems less stringent. More importantly, Lemma 1 gives a *strong law* of large number for $n^{-1} \sum_{t=1}^{n} z_t^2$ under rather mild assumptions, which is one of the key tools for our asymptotic analysis of APE.

**Lemma 2.** *Assume* $\sup_{-\infty < t < \infty} E|\omega_t|^\alpha < \infty$ *for some* $\alpha > 2$ *and*

$$(2.7) \quad \sum_{j \geq k} |c_j| = O(k^{-1}).$$

*Then,*

$$\log \left( \sum_{j=1}^{n-1} x_j^2 \right) = 2 \log n + o(\log n) \text{ a.s.}$$

*Proof.* First note that $x_t = \sum_{j=1}^{t} \eta_j$. Define $N_t = \theta \sum_{j=1}^{t} \omega_j$, where $\theta = \sum_{j=0}^{\infty} c_j$. Then,

$$(2.8) \quad x_t = N_t - S_t,$$

where $S_t = \sum_{j=0}^{t-1} f_j \omega_{t-j}$ with $f_j = \sum_{l=j+1}^{\infty} c_l$. In view of (2.8),

$$(2.9) \quad \sum_{j=1}^{n-1} x_j^2 = \sum_{j=1}^{n-1} N_j^2 - 2 \sum_{j=1}^{n-1} N_j S_j + \sum_{j=1}^{n-1} S_j^2.$$

Since $|f_j| = O(j^{-1})$, Lemma 1 yields

$$\sum_{j=1}^{n-1} S_j^2 = O(n) \text{ a.s.} \tag{2.10}$$

By the law of the iterated logarithm,

$$\sum_{j=1}^{n-1} N_j^2 = O(n^2 \log \log n) \text{ a.s.} \tag{2.11}$$

By Lai and Wei ([12], (3.23)),

$$\liminf_{n \to \infty} \frac{\log \log n}{n^2} \sum_{j=1}^{n-1} N_j^2 > 0 \text{ a.s.} \tag{2.12}$$

Now, Lemma 2 follows directly from (2.9)-(2.12). □

**Remark 2.** By assuming

$$\sum_{j=0}^{\infty} j|c_j| < \infty, \tag{2.13}$$

Proposition 17.3 of Hamilton (1994) gives the limiting distribution of $n^{-2} \sum_{j=1}^{n-1} x_j^2$, which is $\lambda^2 \int_0^1 w(r)^2 dr$, where $\lambda = \sigma_\omega \sum_{j=0}^{\infty} c_j$ and $w(r)$ denotes the standard Brownian motion. This result immediately implies

$$\log\left(\sum_{j=1}^{n-1} x_j^2\right) = 2 \log n + O_p(1). \tag{2.14}$$

Lemma 2 and (2.14) provide different estimates for the difference between $2 \log n$ and $\log(\sum_{j=1}^{n-1} x_j^2)$, but neither is more informative than the other. On the other hand, we have found that the assumption on the coefficients used in Lemma 2, (2.7), seems to be weaker than the one imposed by Hamilton, (2.13). This can be seen by observing that (2.7) is marginally satisfied by $C_1 j^{-2} \leq |c_j| \leq C_2 j^{-2}, C_2 \geq C_1 > 0$, whereas (2.13) is not.

We are now ready to prove the main result of this section.

**Theorem 1.** *Assume that models (1.1), (1.2), and the assumptions of Lemma 2 hold. Also assume that $\sup_{-\infty < t < \infty} E|\varepsilon_t|^{\alpha_0} < \infty$ for some $\alpha_0 > 2$. Then,*

$$\sum_{i=2}^{n} x_{i-1}^2 (\hat{\beta}_{i-1} - \beta)^2 = 2\sigma^2 \log n + o(\log n) \text{ a.s.,} \tag{2.15}$$

*and*

$$\sum_{i=2}^{n} (y_i - \hat{y}_i)^2 = \sum_{i=2}^{n} \varepsilon_i^2 + 2\sigma^2 \log n + o(\log n) \text{ a.s.} \tag{2.16}$$

*Proof.* First note that (2.9)-(2.12) yield

$$\text{(2.17)} \qquad \limsup_{n\to\infty} \frac{n^2}{(\log\log n)\sum_{j=1}^{n} x_j^2} < \infty \text{ a.s.}$$

By Wei ([15], Lemma 2) and (2.7),

$$\text{(2.18)} \qquad E\left|\frac{S_n}{n^{1/2}}\right|^\alpha \leq C_\alpha n^{-\alpha/2}\left(\sum_{j=0}^{n-1} f_j^2\right)^\alpha \leq C_\alpha^* n^{-\alpha/2},$$

where $C_\alpha$ and $C_\alpha^*$ depend only on $\alpha$. (2.18) and the Borel-Cantelli lemma give

$$\text{(2.19)} \qquad S_n = o(n^{1/2}) \text{ a.s.}$$

Since the law of the iterated logarithm implies

$$N_n = O((n\log\log n)^{1/2}) \text{ a.s.},$$

this, (2.8), (2.17), and (2.19) yield

$$\text{(2.20)} \qquad \frac{x_n^2}{\sum_{j=1}^{n} x_j^2} = o(1) \text{ a.s.}$$

In view of (2.20) and Wei ([15], Theorem 3), we have

$$\text{(2.21)} \qquad \sum_{i=2}^{n} x_{i-1}^2(\hat{\beta}_{i-1} - \beta)^2 = \sigma^2 \log\left(\sum_{j=1}^{n-1} x_j^2\right) + o\left(\log\left(\sum_{j=1}^{n-1} x_j^2\right)\right) \text{ a.s.,}$$

As a result, (2.15) follows from Lemma 2 and (2.21); and (2.16) is an immediate consequence of (2.15) and (1.6). □

## 3. An asymptotic expression for the FPE

Assume that models (1.1) and (1.2) hold, $E(\varepsilon_t \omega_t) = \pi$ is a constant independent of $t$, $\sup_{-\infty < t < \infty} E|\varepsilon_t|^{\alpha_0} < \infty$, $\alpha_0 > 2$, and $\sup_{-\infty < t < \infty} E|\omega_t|^\alpha < \infty$, $\alpha > 2$. Then, by the functional central limit theorem, continuous mapping theorem, Ito's formula, and some algebraic manipulations, it can be shown that

$$\text{(3.1)} \qquad \left\{\frac{(\frac{1}{\sqrt{n}}x_n)(\frac{1}{n}\sum_{i=1}^{n-1} x_i \varepsilon_{i+1})}{\frac{1}{n^2}\sum_{i=1}^{n-1} x_i^2}\right\}^2 \\ \Longrightarrow \frac{w_a^2(1)\left(\rho\sigma_\omega \int_0^1 w_a(t)dw_a(t) + \sigma_\theta \int_0^1 w_a(t)dw_b(t)\right)^2}{\left(\int_0^1 w_a^2(t)dt\right)^2},$$

where "$\Longrightarrow$" denotes weak convergence, $(w_a(t), w_b(t))$ is a standard Brownian motion of dimension 2, $\rho = \pi/\sigma_\omega^2$, and $\sigma_\theta^2 = \sigma^2 - \rho^2\sigma_\omega^2$. If we can further show that for some $q > 2$,

$$\text{(3.2)} \qquad E\left|\frac{(\frac{1}{\sqrt{n}}x_n)(\frac{1}{n}\sum_{i=1}^{n-1} x_i\varepsilon_{i+1})}{\frac{1}{n^2}\sum_{i=1}^{n-1} x_i^2}\right|^q = O(1),$$

then, in view of (3.1), (3.2), and (1.8),

$$
(3.3) \quad nE\{x_n^2(\hat{\beta}_n - \beta)^2\} = E\left\{\frac{w_a^2(1)\left(\rho\sigma_\omega \int_0^1 w_a(t)dw_a(t) + \sigma_\theta \int_0^1 w_a(t)dw_b(t)\right)^2}{\left(\int_0^1 w_a^2(t)dt\right)^2}\right\} + o(1).
$$

In the rest of this section, we provide sufficient conditions to ensure (3.2). In addition, the expectation on the rhs of (3.3) is investigated (Corollary 1). Let us start with a useful lemma.

**Lemma 3.** *Let $F_{t,m,\mathbf{a}_m}(.)$ be the distribution function of $\sum_{j=1}^m a_j\omega_{t+1-j}$, where $\mathbf{a}_m = (a_1, \ldots, a_m)'$. There are some positive numbers $\kappa$, $\iota$, and $M$ such that for all $m \geq 1, -\infty < t < \infty$ and $\|\mathbf{a}_m\|^2 = \sum_{j=1}^m a_j^2 = 1$,*

$$
(3.4) \quad |F_{t,m,\mathbf{a}_m}(x) - F_{t,m,\mathbf{a}_m}(y)| \leq M |x-y|^\kappa,
$$

*as $|x-y| \leq \iota$. Then, for any $q > 0$,*

$$
(3.5) \quad E\left\{\left(\frac{1}{n^2}\sum_{j=1}^{n-1} x_j^2\right)^{-q}\right\} = O(1).
$$

*Proof.* The proof is closely related to the one given in Ing ([9], Lemma 1), with the assumption there being strengthened to (3.4). First note that

$$
(3.6) \quad \frac{1}{n^2}\sum_{i=1}^{n-1} x_i^2 \geq \frac{1}{n^2}\sum_{i=n\delta}^{n-1} x_i^2 = \frac{\delta}{n}\sum_{i=n\delta}^{n-1} \frac{x_i^2}{n\delta} \geq \frac{\delta}{n}\sum_{i=n\delta}^{n-1} \frac{x_i^2}{i},
$$

where $0 < \delta < 1$, and without loss of generality, $n\delta$ is assumed to be a positive integer. Rearranging the series on the rhs of (3.6), one obtains

$$
(3.7) \quad \frac{\delta}{n}\sum_{j=0}^{\frac{(1-\delta)n}{lq}-1}\sum_{i=0}^{lq-1}\frac{x^2_{n\delta+\frac{(1-\delta)n}{lq}i+j}}{n\delta + \frac{(1-\delta)n}{lq}i + j},
$$

where $l > \max[2/\kappa, 1/q, (1/q)\{(1/\delta)-1)\}]$ and for simplifying the discussion, $lq$ and $\{(1-\delta)n\}/(lq)$ are also assumed to be positive integers. By the convexity of function $x^{-q}, x > 0$,

$$
(3.8) \quad \left(\frac{1}{n^2}\sum_{i=1}^{n-1} x_i^2\right)^{-q} \leq \left\{\frac{(1-\delta)\delta}{lq}\right\}^{-q}\frac{lq}{(1-\delta)n}
$$
$$
\times \sum_{j=0}^{\frac{(1-\delta)n}{lq}-1}\left\{\sum_{i=0}^{lq-1}\frac{x^2_{n\delta+\frac{(1-\delta)n}{lq}i+j}}{n\delta + \frac{(1-\delta)n}{lq}i + j}\right\}^{-q}.
$$

In view of (3.8), if one can show that for some positive number $C$ independent of $j$, the following inequality,

$$
(3.9) \quad E\left\{\sum_{i=0}^{lq-1}\frac{x^2_{n\delta+\frac{(1-\delta)n}{lq}i+j}}{n\delta + \frac{(1-\delta)n}{lq}i + j}\right\}^{-q} \leq C < \infty,
$$

holds for all $j = 0, 1, \ldots, \{(1-\delta)n/(lq)\} - 1$ as $n$ is large enough, then (3.5) follows. The rest of the proof only focuses on the case where $j = 0$, because the same argument can be easily applied to other $j$'s.

For $i = 0, \ldots, lq - 1$, define

$$(3.10) \qquad Y_{n,i} = \left\{ n\delta + \frac{(1-\delta)n}{lq} i \right\}^{-1/2} x_{n\delta + \frac{(1-\delta)n}{lq} i},$$

$$(3.11) \qquad W_{n,i} = \left\{ n\delta + \frac{(1-\delta)n}{lq} i \right\}^{-1/2} \sum_{m=0}^{\frac{(1-\delta)n}{lq} - 1} \bar{f}_m \omega_{n\delta + \frac{(1-\delta)n}{lq} i - m},$$

where $\bar{f}_j = \sum_{l=0}^{j} c_l$, and

$$(3.12) \qquad F_{n,i} = Y_{n,i} - W_{n,i}.$$

(Note that $x_t = \sum_{j=0}^{t-1} \bar{f}_j \omega_{t-j}$.) Then,

$$(3.13) \quad \begin{aligned} E\left(\sum_{i=0}^{lq-1} Y_{n,i}^2\right)^{-q} &= \int_0^\infty Pr\left\{ \left(\sum_{i=0}^{lq-1} Y_{n,i}^2\right)^{-q} > t \right\} dt \\ &= \int_0^\infty Pr\left(\sum_{i=0}^{lq-1} Y_{n,i}^2 < t^{-1/q}\right) dt \\ &\leq \int_0^\infty Pr\left(-t^{-1/(2q)} < Y_{n,i} < t^{-1/(2q)}, \quad i = 0, \ldots, lq-1\right) dt \\ &= \int_0^\infty E\left\{ E\left(\prod_{i=0}^{lq-1} I_{A_{n,i}} \middle| F_{n,lq-1}, W_{n,i}, F_{n,i}, i = 0, \ldots, lq-2\right) \right\} dt, \end{aligned}$$

where $A_{n,i} = \{-t^{-1/(2q)} < Y_{n,i} < t^{-1/(2q)}\}$. In view of (3.10)-(3.12), for $0 \leq p \leq lq - 1$, $0 \leq i \leq p$, and $0 \leq j \leq p - 1$, $W_{n,p}$ is independent of $(F_{n,i}, W_{n,j})$. In addition, $\text{var}(W_{n,i}) > \zeta > 0$, where $i = 0, \ldots, lq - 1$ and $\zeta$ is a positive number independent of $n$ and $i$. According to these facts, (3.4), and arguments similar to those used in (3.10) and (3.11) of Ing [9], there exist some positive numbers $0 < C' < \infty, 0 < s < \infty$, and a positive integer $N_0$ such that for all $n \geq N_0$ and all $t \geq s$,

$$(3.14) \qquad E\left(\prod_{i=0}^{lq-1} I_{A_{n,i}}\right) \leq C' t^{-(\kappa l)/2}.$$

Since, by construction, $l > 2/\kappa$, (3.13) and (3.14) guarantee that for $n > N_0$,

$$E\left(\sum_{i=0}^{lq-1} Y_{n,i}^2\right)^{-q} \leq s + C' \int_s^\infty t^{-(\kappa l)/2} dt < \infty,$$

which yields (3.9). □

Lemma 4 below shows that (3.4) is easily found in many time series applications.

**Lemma 4.** *If $\omega_t$'s are i.i.d. random variables satisfying $E(\omega_1) = 0, E(\omega_1^2) = \sigma_\omega^2 > 0$, and $E(|\omega_1|^\alpha) < \infty$ for some $\alpha > 2$. Assume also that for some positive constant $M_0 < \infty$,*

$$(3.15) \qquad \int_{-\infty}^\infty |\varphi(t)| dt \leq M_0,$$

where $\varphi(t) = E(e^{it\omega_1})$ is the characteristic function of $\omega_1$. Then, for all $-\infty < t < \infty$, $m \geq 1$ and $\|\mathbf{a}_m\| = 1$, there is a finite positive constant $M_1$ such that

$$\sup_{-\infty < x < \infty} f_{t,m,\mathbf{a}_m}(x) < M_1, \tag{3.16}$$

where $f_{t,m,\mathbf{a}_m}(\cdot)$ is the density function of $\sum_{j=1}^m a_j \omega_{t+1-j}$. As a result, (3.4) follows.

*Proof.* The proof is inspired by the ideas of Feller ([7], p. 516), which deal with the special case, $a_j = m^{-1/2}$ for all $j = 1, \ldots, m$. Without loss of generality, assume $\sigma_\omega^2 = 1$. Denote $Y = \sum_{j=1}^m a_j \omega_{t+1-j}$. Then, $\varphi_Y(t) = E(e^{itY}) = \prod_{j=1}^m \varphi_j(a_j t)$. By Chow and Teicher ([6], Theorem 8.4.1),

$$\varphi(a_j t) = 1 - \frac{a_j^2 t^2}{2} + o(a_j^2 t^2),$$

as $a_j^2 t^2 \to 0$. This gives for $|a_j t| < \delta_1^*$, where $\delta_1^*$ is some small positive constant,

$$|\varphi(a_j t)| \leq 1 - \frac{a_j^2 t^2}{4}. \tag{3.17}$$

On the other hand, since (3.15) yields $|\varphi(t)| \to 0$ as $|t| \to \infty$, by Chow and Teicher ([6], Corollary 8.4.2), $|\varphi(t)| < 1$ for all $t \neq 0$, and hence for all $|t| \geq \delta_1^*$ (with $\delta_1^*$ defined above),

$$|\varphi(t)| < \theta_1, \tag{3.18}$$

where $\theta_1$ is some positive constant $< 1$. Now, by (3.17),

$$\begin{aligned}\int_{-\infty}^{\infty} \prod_{j=1}^m |\varphi(a_j t)| dt &\leq \int_{|t| < \frac{\delta_1^*}{O_m}} e^{\frac{-t^2}{4}} dt + \int_{|t| \geq \frac{\delta_1^*}{O_m}} \prod_{j=1}^m |\varphi(a_j t)| dt \\ &\leq \int_{-\infty}^{\infty} e^{\frac{-t^2}{4}} dt + \int_{|t| \geq \frac{\delta_1^*}{O_m}} \prod_{j=1}^m |\varphi(a_j t)| dt,\end{aligned} \tag{3.19}$$

where $O_j$ is a permutation of $|a_j|$ satisfying $O_m \geq O_{m-1} \geq \cdots \geq O_1$. For $t \geq \delta_1^*/O_m$, (3.17), (3.18) and the fact that

$$\theta_1 = 1 - (1 - \theta_1) \leq 1 - \frac{4(1-\theta_1)}{\delta_1^{*2}} \frac{a_j^2 \delta_1^{*2}}{4 O_m^2}$$

imply

$$|\varphi(a_j t)| \leq \max\{1 - \frac{a_j^2 t^2}{4}, \theta_1\} \leq \max\{1 - \frac{a_j^2 \delta_1^{*2}}{4 O_m^2}, \theta_1\} \leq 1 - \xi \frac{a_j^2 \delta_1^{*2}}{4 O_m^2}, \tag{3.20}$$

where $0 < \xi < \min\{1, 4(1-\theta_1)/\delta_1^{*2}\}$. In view of (3.20) and the fact that $\sum_{j=1}^{m-1} O_j^2 = 1 - O_m^2$,

$$\begin{aligned}\int_{|t| \geq \frac{\delta_1^*}{O_m}} \prod_{j=1}^m |\varphi(a_j t)| dt &\leq \frac{1}{O_m} \int_{-\infty}^{\infty} e^{-\frac{\xi \delta_1^{*2}}{4 O_m^2} \sum_{j=1}^{m-1} O_j^2} |\varphi(t)| dt \\ &= e^{\frac{\xi \delta_1^{*2}}{4}} \frac{1}{O_m} e^{\frac{-\xi \delta_1^{*2}}{4 O_m^2}} \int_{-\infty}^{\infty} |\varphi(t)| dt \\ &\leq e^{\frac{\xi \delta_1^{*2}}{4}} \sup_{x \geq 1} x e^{\frac{-\xi \delta_1^{*2} x^2}{4}} M_0 < \infty.\end{aligned} \tag{3.21}$$

By (3.21), (3.19), and the fact that

$$\sup_{-\infty<x<\infty} f_{t,m,\mathbf{a}_m}(x) \leq \frac{1}{2\pi}\int_{-\infty}^{\infty}\prod_{j=1}^{m}|\varphi(a_jt)|dt,$$

(3.16) follows. In addition, it is not difficult to see that (3.4) can be deduced from (3.16). □

In the following lemma, some moment bounds for $(1/\sqrt{n})x_n$ and $(1/n)\sum_{i=1}^{n-1} x_i \times \varepsilon_{i+1}$, are obtained.

**Lemma 5.** *Assume models (1.1) and (1.2), with $\sup_t E(|\varepsilon_t|^q) < \infty$ and $\sup_t E(|\omega_t|^q) < \infty$, for some $q \geq 2$. Then,*

(3.22)      (i)      $\displaystyle\sup_{n\geq 1} E\left(\left|\frac{1}{\sqrt{n}}x_n\right|^q\right) < \infty,$

(3.23)      (ii)      $\displaystyle\sup_{n\geq 1} E\left(\left|\frac{1}{n}\sum_{i=1}^{n-1} x_i\varepsilon_{i+1}\right|^q\right) < \infty.$

*Proof.* The proof of Lemma 5 is similar to that of Ing ([9], Lemma 1). The details are omitted. □

Armed with the previous results, (3.2) is proved in the following theorem.

**Theorem 2.** *Assume that (1.1), (1.2), (3.4), $\sup_t E(|\varepsilon_t|^q) < \infty$, and $\sup_t E(|\omega_t|^q) < \infty$ are satisfied, where $q > 4$. Then, (3.2) holds. If we further assume that $E(\varepsilon_t\omega_t) = \pi$ is a constant independent of $t$, then (3.3) follows.*

*Proof.* By Lemmas 3 and 5, (3.1), and an argument similar to the one used in [9], Theorem 1, the claimed results can be obtained. □

The FPE of the least squares predictor is obtained in Corollary 1 below.

**Corollary 1.** *Assume that (2.7) and all assumptions of Theorem 2 hold. Then, (1.9) follows.*

*Proof.* By (2.15), (3.2), and Minkowski's inequality,

$$(3.24) \qquad \lim_{n\to\infty} \frac{1}{\log n}\sum_{i=m^*}^{n} E\{x_{i-1}^2(\hat{\beta}_{i-1}-\beta)^2\} = 2\sigma^2,$$

where $m^*$ is some positive integer independent of $n$. Now, (1.9) is guaranteed by (3.3) and (3.24). □

Corollary 1 and Theorem 1 together indicate an interesting result that the term of order $\log n$ in the APE and the term of order $n^{-1}$ in the FPE share the same constant, $2\sigma^2$. For applications of this type of results to model selection problems, see [11]. Corollary 1 also shows that the FPE of the least squares predictor is not affected by the contemporary correlation between $\varepsilon_t$ and $\omega_t$. This is a somewhat unexpected feature because the least squares estimate itself does not possess this property. More specifically, by direct calculations, we have

$$(3.25) \qquad n(\hat{\beta}_n - \beta) \Longrightarrow \frac{1}{\lambda}\frac{\rho\sigma_\omega \int_0^1 w_a(t)dw_a(t) + \sigma_\theta \int_0^1 w_a(t)dw_b(t)}{\int_0^1 w_a^2(t)dt},$$

and

$$(3.26) \quad n^2(\hat{\beta}_n - \beta)^2 \Longrightarrow \frac{1}{\lambda^2} \frac{\left(\rho\sigma_\omega \int_0^1 w_a(t)dw_a(t) + \sigma_\theta \int_0^1 w_a(t)dw_b(t)\right)^2}{\left(\int_0^1 w_a^2(t)dt\right)^2},$$

where $\lambda$ is defined in Remark 2. By (3.26), an argument similar to that used in the proof of Theorem 2, and some algebraic manipulations,

$$(3.27) \quad \begin{aligned} \lim_{n\to\infty} n^2 E(\hat{\beta}_n - \beta)^2 &= E\left\{\frac{1}{\lambda^2} \frac{\left(\rho\sigma_\omega \int_0^1 w_a(t)dw_a(t) + \sigma_\theta \int_0^1 w_a(t)dw_b(t)\right)^2}{\left(\int_0^1 w_a^2(t)dt\right)^2}\right\} \\ &= \frac{\rho^2}{\iota^2} E\left(\frac{\int_0^1 w_a(t)dw_a(t)}{\int_0^1 w_a^2(t)dt}\right)^2 + \frac{\sigma_\theta^2}{\iota^2 \sigma_\omega^2} E\left(\frac{1}{\int_0^1 w_a^2(t)dt}\right), \end{aligned}$$

where $\iota^2 = \lambda^2 \sigma_\omega^{-2}$. Ing ([9], (4.3)) showed that

$$(3.28) \quad E\left(\frac{\int_0^1 w_a(t)dw_a(t)}{\int_0^1 w_a^2(t)dt}\right)^2 \doteq 13.3.$$

By (3.6.4) and (3.6.5) of Arató and using a numerical integration method,

$$(3.29) \quad E\left(\frac{1}{\int_0^1 w_a^2(t)dt}\right) \doteq 5.6.$$

Consequently, (3.27)-(3.29) imply

$$(3.30) \quad \lim_{n\to\infty} n^2 E(\hat{\beta}_n - \beta)^2 \doteq \frac{\rho^2}{\iota^2} 13.3 + \frac{\sigma_\theta^2}{\iota^2 \sigma_\omega^2} 5.6,$$

which obviously varies with the strength of dependence between $\varepsilon_t$ and $\omega_t$. In particular, if $\sigma^2 = \sigma_\omega^2$, then $\rho = \mathrm{corr}(\varepsilon_t, \omega_t)$ and (3.30) can be rewritten as

$$(3.31) \quad \lim_{n\to\infty} n^2 E(\hat{\beta}_n - \beta)^2 \doteq \frac{1}{\iota^2}[\rho^2 13.3 + (1-\rho^2)5.6].$$

As observed in (3.31), the larger the magnitude of the correlation between $\varepsilon_t$ and $\omega_t$ is, the larger the mean squared error of the least squares estimate is, a result new to the literature.

As a final remark, we note that the square of the normalized estimate, $n^2(\hat{\beta}_n - \beta)^2$, and the square of normalized regressor, $x_n^2/n$, are not asymptotically uncorrelated. To see this, observe that $\lim_{n\to\infty} E(x_n^2/n) = \lambda^2$, which together with (3.30) and Corollary 1, gives

$$\begin{aligned} \lim_{n\to\infty} E\left(\frac{x_n^2}{n}\right) E\left\{n^2(\hat{\beta}_n - \beta)^2\right\} &\doteq 13.3\rho^2\sigma_\omega^2 + 5.6\sigma_\theta^2 \\ &= 5.6\sigma^2 + 7.7\rho^2\sigma_\omega^2 > 2\sigma^2 \\ &= \lim_{n\to\infty} E\left\{\frac{x_n^2}{n} n^2(\hat{\beta}_n - \beta)^2\right\}. \end{aligned}$$

Therefore, $x_n^2/n$ and $n^2(\hat{\beta}_n - \beta)^2$ are (asymptotically) negatively correlated, which suggests that larger variation of $x_n$ can yield a better estimation result. It is worth mentioning that this special feature does not exist for the (asymptotically) stationary regressor. For example, when $x_t = \varsigma x_{t-1} + \eta_t$, with $|\varsigma| < 1$, following an argument used in Ing [10], it can be shown that

$$\lim_{n\to\infty} E(x_n^2) E\left\{[\sqrt{n}(\hat{\beta}_n - \beta)]^2\right\} = \lim_{n\to\infty} E\left\{x_n^2 n(\hat{\beta}_n - \beta)^2\right\} = \sigma^2.$$

Therefore, the square of the normalized estimate, $n(\hat{\beta}_n - \beta)^2$, and the square of the (normalized) regressor, $x_n^2$, are asymptotically uncorrelated in this case.